\documentclass[a4paper,11pt]{amsart}

\usepackage{verbatim}
\usepackage{amsmath}
\usepackage{amssymb}
\usepackage{amsfonts}
\usepackage{amsthm}
\usepackage{mathrsfs}
\usepackage{tikz}
\usetikzlibrary{matrix,arrows}
\usepackage[all]{xy}
\usepackage{mathrsfs}
\usepackage{a4wide}
\usepackage{enumerate}
\usepackage{hyperref}

\theoremstyle{plain}

\newtheorem{theorem}{Theorem}
\newtheorem{lemma}[theorem]{Lemma}
\newtheorem{proposition}[theorem]{Proposition}
\newtheorem{corollary}[theorem]{Corollary}

\newtheorem*{conjecture1}{Conjecture}
\newtheorem{rem}[theorem]{Remark}

\theoremstyle{definition}
\newtheorem{definition}[theorem]{Definition}
\newtheorem{example}[theorem]{Example}

\theoremstyle{remark}
\newtheorem{remark}[theorem]{Remark}

\newcommand{\G}{\mathcal{G}}

\newcommand{\ucong}{\rotatebox{90}{$\cong$}}

\title{Fibred coarse embeddings, a-T-menability and the coarse analogue of the Novikov conjecture} 
\date{September 2014}
\author{Martin Finn-Sell}
\email[M. Finn-Sell]{ mfinnse@uni-math.gwdg.de}
\curraddr{Mathematische Institut, Georg-August-Universit\"{a}t G\"{o}ttingen, 3-5 Bunsen Stra{\ss}e, 37073, G\"{o}ttingen, Deutschland}

\begin{document}
\bibliographystyle{plain}

\begin{abstract}
The connection between the coarse geometry of metric spaces and analytic properties of topological groupoids is well known. One of the main results of Skandalis, Tu and Yu is that a space admits a coarse embedding into Hilbert space if and only if a certain associated topological groupoid is a-T-menable. This groupoid characterisation then reduces the proof that the coarse Baum-Connes conjecture holds for a coarsely embeddable space to known results for a-T-menable groupoids. The property of admitting a fibred coarse embedding into Hilbert space was introduced by Chen, Wang and Yu to provide a property that is sufficient for the maximal analogue to the coarse Baum-Connes conjecture and in this paper we connect this property to the traditional coarse Baum-Connes conjecture using a restriction of the coarse groupoid and homological algebra. Additionally we use this results to give a characterisation of the a-T-menability for residually finite discrete groups.
\end{abstract}

\maketitle


\section{Introduction}

The application of coarse methods to algebraic topological problems is well known \cite{MR866507,MR1728880}. One possible method utilised for such problems is a \textit{higher index} theorem that allows us to calculate refined large scale information from small scale topological or analytic data. This process can be succinctly phrased using the language of K-theory and K-homology and can be encoded into the coarse geometric version of the well-known Baum-Connes conjecture, which asks whether or not a certain assembly map:
\begin{equation*}
\mu_{*}:KX_{*}(X)=\lim_{R>0}K_{*}(P_{R}(X)) \longrightarrow K_{*}(C^{*}(X))
\end{equation*}
is an isomorphism for all uniformly discrete spaces $X$ with bounded geometry. One approach to this conjecture for suitable metric spaces $X$ is via the concept of a \textit{coarse embedding into Hilbert space}. A seminal paper by Yu \cite{MR1728880} first showed the importance of coarse embeddings into Hilbert space by proving that this is a sufficient condition for the coarse Baum-Connes assembly map $\mu_{*}$ to be an isomorphism.

In this paper we study the relationship between \textit{fibred coarse embedding into Hilbert space}, first introduced by Chen, Wang and Yu \cite{MR3116568}, and the coarse analogue of the strong Novikov conjecture. Intuitively, a space admits a fibred coarse embedding into Hilbert space if for each scale it is acceptable to forget bounded portions of the space and embed what remains locally and compatibly into Hilbert space. This is made precise in Definition \ref{def:FCE}. This property was used in \cite{MR3116568} to prove a maximal analogue of the work of \cite{MR1728880}, that is the existence of a fibred coarse embedding into Hilbert space implies that the maximal coarse Baum-Connes assembly map is an isomorphism for any uniformly discrete metric space with bounded geometry.

Another approach to these questions was considered in \cite{mypub1}, in which a conjecture known as the \textit{boundary coarse Baum-Connes conjecture} was defined for uniformly discrete bounded geometry spaces. This conjecture, defined via groupoids, was designed to capture the structure of a space at infinity. We explain the relationship between this conjecture and the work of \cite{MR3116568} and give a new method to prove the maximal coarse Baum-Connes conjecture in this instance that has additional consequences. We explore these consequences in Section \ref{sect:apps}.

The strategy of this paper is to convert coarse assembly problems into groupoid assembly problems was pioneered in \cite{MR1905840}:

\begin{theorem}\label{Thm:MT1}
Let $X$ be a uniformly discrete space with bounded geometry that admits a fibred coarse embedding into Hilbert space. Then the associated boundary groupoid $G(X)|_{\partial\beta X}$ is a-T-menable.
\end{theorem}

We define the boundary groupoid $G(X)|_{\partial\beta X}$ in Definition \ref{def:bdry}.

This result gives us access to the tools developed in \cite{mypub1} concerning the boundary coarse Baum-Connes conjecture and in in Section \ref{sect:apps} we provide applications of Theorem \ref{Thm:MT1} to the coarse analogue of the strong Novikov conjecture for uniformly discrete spaces with bounded geometry (Theorem \ref{thm:mcor1}), as well as considering the situation concerning ghost operators within Roe algebras associated to coarsely embeddable spaces (Corollary \ref{thm:IT3}). Lastly we give a characterisation of a-T-menability for residually finite groups via box spaces (Theorem \ref{thm:cor2}).

\section*{Acknowledgements}
The author would like to thank his supervisor Nick Wright and J\'{a}n \v{S}pakula for some helpful conversations and comments. I would also like to thank the referee for pointing out numerous errors in the previous versions as well as their patience during the correction process and finally Rufus Willett for sharing the idea that motivates Section \ref{sect:whathappens}.

\section{Fibred coarse embeddings and groupoids associated to coarse spaces}\label{sect:coarse}
We recap now the important details surrounding Theorem 5.4 from \cite{MR1905840} and the notion of a fibred coarse embedding \cite{MR3116568}.
\subsection{Coarse and fibred coarse embeddings}
We first make precise the metric spaces that we will study in this paper.

\begin{definition}
Let $X$ be a metric space. Then $X$ is said to be \textit{uniformly discrete} if there exists $c>0$ such that for every pair of distinct points $x\not = y \in X$ the distance $d(x,y)>c$. Additionally, $X$ is said have \textit{bounded geometry} if for every $R>0$ there exists $N_{R}>0$ such that for every $x \in X$ the cardinality of the ball of radius $R$ about $x$ is smaller than $N_{R}$.
\end{definition}

And now recall the concept of a coarse embedding, to motivate the definition of a fibred coarse embedding.

\begin{definition}
A metric space $X$ is said to admit a coarse embedding into Hilbert space $\mathcal{H}$ if there exist maps $f:X \rightarrow \mathcal{H}$,  and non-decreasing $\rho_{1},\rho_{2}:\mathbb{R}_{+} \rightarrow \mathbb{R}$ such that:
\begin{enumerate}
\item for every $x,y \in X$, $\rho_{1}(d(x,y)) \leq \Vert f(x) - f(y) \Vert \leq \rho_{2}(d(x,y))$;
\item for each $i$, we have $\lim_{r \rightarrow \infty}\rho_{i}(r) = +\infty$.
\end{enumerate}
\end{definition}

One application of admitting a coarse embedding into Hilbert space is the main result of \cite{MR1728880, MR1905840}.

\begin{theorem}
Let $X$ be a uniformly discrete metric space with bounded geometry that admits a coarse embedding into Hilbert space. Then the coarse Baum-Connes conjecture holds for $X$, that is the assembly map $\mu_{*}$ is an isomorphism.\qed
\end{theorem}

Many metric spaces admit coarse embeddings; the property is a suitably flexible one as it is implied by many other coarse properties such as finite asymptotic dimension \cite[Example 11.5]{MR2007488}, or the weaker properties of finite decomposition complexity \cite{MR2947546} and property A \cite{MR1728880}. The primary application is to finitely generated discrete groups via the descent principle that is originally due to Higson \cite{MR1779613}, then adapted by Skandalis, Tu and Yu in \cite{MR1905840}:

\begin{corollary}
Let $\Gamma$ be a finitely generated discrete group admitting a coarse embedding into Hilbert space. Then the (strong) Novikov conjecture holds for $\Gamma$.\qed
\end{corollary}

For more information on the strong Novikov conjecture and its connections to coarse geometry see \cite{MR1388300, MR1905840}.

There are spaces that do not admit a coarse embedding into Hilbert space, such as expander graphs \cite{MR2569682}. Certain types of expander graph are known to be counterexamples to the Baum-Connes conjecture \cite{higsonpreprint,MR1911663,explg1,explg2,MR2568691} and within this class, it is possible to prove some positive results under hypotheses on the structure of the expander; it is known that certain families do satisfy the coarse analogue of the Novikov conjecture and the maximal coarse Baum-Connes conjecture \cite{MR2568691,explg1}. The notion of a fibred coarse embedding was introduced in \cite{MR3116568} where the authors gave results that partially explained these phenomena. The following is Definition 2.1 from \cite{MR3116568}.

\begin{definition}\label{def:FCE}
A metric space $X$ is said to admit a \textit{fibred coarse embedding} into Hilbert space if there exists
\begin{itemize}
\item a field of Hilbert spaces $\lbrace H_{x} \rbrace_{x \in X}$ over X;
\item a section $s: X \rightarrow \sqcup_{x \in X}H_{x}$ (i.e $s(x) \in H_{x}$);
\item two non-decreasing functions $\rho_{1}, \rho_{2}$ from $[0,\infty)$ to $(-\infty, +\infty)$ such that $\lim_{r\rightarrow \infty}\rho_{i}(r)=\infty$ for $i=1,2$;
\item A reference Hilbert space $H$
\end{itemize}
such that for any $r>0$ there exists a bounded subset $K_{r}\subset X$ and a trivalisation:
\begin{equation*}
t_{C}: \sqcup_{x\in C} H_{x} \rightarrow C\times H
\end{equation*}
for each $C \subset X \setminus K_{r}$ of diameter less than $r$. We ask that this map $t_{C}|_{H_{x}}=t_{C}(x)$ is an affine isometry from  $H_{x}$ to $H$, satisfying:
\begin{enumerate}
\item for any $x,y \in C$, $\rho_{1}(d(x,y))\leq \Vert t_{C}(x)(s(x)) - t_{C}(y)(s(y)) \Vert \leq \rho_{2}(d(x,y))$;
\item for any two subsets $C_{1},C_{2} \subset X\setminus K_{r}$ of diameter less than $r$ and nonempty intersection $C_{1}\cap C_{2}$ there exists an affine isometry $t_{C_{1},C_{2}}:H \rightarrow H$ such that $t_{C_{1}}(x)t_{C_{2}}(x)^{-1}=t_{C_{1},C_{2}}$ for all $x \in C_{1}\cap C_{2}$.
\end{enumerate}
Let $\mathcal{F}=\lbrace K_{r} \rbrace_{r}$. In this instance we say $X$ fibred coarsely embeds into Hilbert space with respect to $\mathcal{F}$.  
\end{definition}

It is possible to improve this definition in the context of \textit{sequences of finite metric spaces} and the infinite metric spaces defined from them using the following construction:
\begin{definition}\label{def:coarselydisconnected}
Let $\lbrace X_{i} \rbrace_{i \in \mathbb{N}}$ be a sequence of finite metric spaces that are uniformly discrete with bounded geometry uniformly in the index $i$, such that $\vert X_{i} \vert \rightarrow \infty$ in $i$. Then we can form the \textit{coarse disjoint union} $X$ with underlying set $\sqcup X_{i}$, metric $d$ given by the metric on each component and $d(X_{i},X_{j})\rightarrow \infty$ as $i+j\rightarrow \infty$. Any such metric is proper and unique up to coarse equivalence.
\end{definition}

In this case, we can refine Definition \ref{def:FCE} to the following:

\begin{definition}\label{def:FCE2}
The coarse disjoint union $X$ of a sequence of finite metric spaces $\lbrace X_{i} \rbrace$ is said to admit a \textit{fibred coarse embedding} into Hilbert space if there exists
\begin{itemize}
\item a field of Hilbert spaces $\lbrace H_{x} \rbrace_{x \in X}$ over X;
\item a section $s: X \rightarrow \sqcup_{x \in X}H_{x}$ (i.e $s(x) \in H_{x}$);
\item two non-decreasing functions $\rho_{1}, \rho_{2}$ from $[0,\infty)$ to $(-\infty, +\infty)$ such that $\lim_{r\rightarrow \infty}\rho_{i}(r)=\infty$ for $i=1,2$;
\item a non-decreasing sequence of numbers $0\leq l_{0} \leq l_{1} \leq ... \leq l_{i} \leq ...$ with $\lim_{i\rightarrow \infty} l_{i} = \infty$;
\item a reference Hilbert space $H$
\end{itemize}
such that for each $x \in X_{i}$ there exists a trivialisation:
\begin{equation*}
t_{x}: \sqcup_{y \in  B_{l_{i}}(x)} H_{y}  \rightarrow B_{l_{i}}(x)\times H
\end{equation*}
such that the map $t_{x}(z)$ is an affine isometry from  $H_{z}$ to $H$, satisfying:
\begin{enumerate}
\item $\rho_{1}(d(z,z^{'}))\leq \Vert t_{x}(z)(s(z)) - t_{x}(z^{'})(s(z^{'})) \Vert \leq \rho_{2}(d(z,z^{'}))$ for any $z,z^{'} \in B_{l_{i}}(x), x \in X_{i}, i\in \mathbb{N}$;
\item for any $x,y\in X_{i}$ such that $B_{l_{i}}(x) \cap B_{l_{i}}(y)\not = \emptyset$ there exists an affine isometry $t_{xy}:H \rightarrow H$ such that $t_{x}(z)=t_{xy}t_{y}(z)$ for all $z \in B_{l_{i}}(x) \cap B_{l_{i}}(y)$.
\end{enumerate}

\end{definition}

In particular many expanders are known to satisfy this property, such as those coming from a-T-menable discrete groups or those of large girth \cite{MR3116568,MR2568691}.

The main application of this property is the main result of \cite{MR3116568}:

\begin{theorem}\label{Thm:FCEMR}
Let $X$ be a bounded geometry metric space admitting a fibred coarse embedding into Hilbert space. Then the maximal coarse Baum-Connes conjecture holds for $X$.\qed
\end{theorem}

We give a different proof of this in Section \ref{Sect:apps}.

\subsection{Groupoids and coarse properties}\label{sect:coarsegroupoids}
Groupoids play an integral role in the constructions we adapt from \cite{MR1905840}. Below we recap some basic properties of groupoids before giving an outline of the construction of the coarse groupoid associated to a coarse space $X$.

A groupoid $\G$ is a \textit{topological groupoid} if both $\G$ and $\G^{(0)}$ are topological spaces, and the maps $r,s, ^{-1}$ and the composition are all continuous. A Hausdorff, locally compact topological groupoid $\G$ is \textit{proper} if $(r,s)$ is a proper map and \textit{\'etale} or \textit{r-discrete} if the map $r$ is a local homeomorphism. When $\G$ is \'etale, $s$ and the product are also local homeomorphisms, and $\G^{(0)}$ is an open subset of $\G$.

Let $X$ be a uniformly discrete bounded geometry metric space. We want to define a groupoid with the property that it captures the coarse information associated to $X$. To do this effectively we need to define what we mean by a \textit{coarse structure} that is associated to a metric. The details of this can be found in \cite{MR2007488}.

\begin{definition}
Let $X$ be a set and let $\mathcal{E}$ be a collection of subsets of $X \times X$. If $\mathcal{E}$ has the following properties:
\begin{enumerate}
\item $\mathcal{E}$ is closed under finite unions;
\item $\mathcal{E}$ is closed under taking subsets;
\item $\mathcal{E}$ is closed under the induced product and inverse that comes from the pair groupoid product on $X \times X$;
\item $\mathcal{E}$ contains the diagonal.
\end{enumerate}
Then we say $\mathcal{E}$ is a \textit{coarse structure} on $X$ and we call the elements of $\mathcal{E}$ \textit{entourages}. If in addition $\mathcal{E}$ contains all finite subsets then we say that $\mathcal{E}$ is \textit{weakly connected}.
\end{definition}

\begin{example}\label{ex:MCS}
Let $X$ be a metric space. Then consider the collection $\mathcal{S}$ of the $R$-neighbourhoods of the diagonal in $X\times X$; that is, for every $R>0$ the set:
\begin{equation*}
\Delta_{R}=\lbrace (x,y) \in X \times X | d(x,y)\leq R \rbrace
\end{equation*}
Let $\mathcal{E}$ be the coarse structure generated by $\mathcal{S}$. This is called the \textit{metric coarse structure} on $X$. If  $X$ is a uniformly discrete bounded geometry  metric space this coarse structure is uniformly locally finite, proper and weakly connected \cite{MR1905840}.
\end{example}

Let $X$ be a uniformly discrete metric space with bounded geometry. We denote by $\beta X$ the Stone-\v{C}ech compactification of $X$ (similarly with $\beta (X\times X)$).

Define $G(X):=\bigcup_{R>0}\overline{\Delta_{R}} \subseteq \beta(X\times X)$. Then $G(X)$ is a locally compact, Hausdorff topological space. To equip it with a product and inverse we would ideally consider the natural extension of the pair groupoid product on $\beta X \times \beta X$.  We remark that the map $(r,s)$ from $X \times X$ extends first to an inclusion into $\beta X \times \beta X$ and universally to $\beta( X \times X)$, giving a map $\overline{(r,s)}: \beta (X\times X) \rightarrow \beta X \times \beta X$. We can restrict this map to each entourage $E \in \mathcal{E}$ allowing us to map the set $G(X)$ to $\beta X \times \beta X$. The following is Corollary 10.18 \cite{MR2007488}

\begin{lemma}\label{Lem:CorRoe}
Let $X$ be a uniformly discrete bounded geometry metric space, let $E$ be any entourage and let $\overline{E}$ its closure in $\beta(X \times X)$. Then the inclusion $E \rightarrow X \times X$ extends to a topological embedding $\overline{E} \rightarrow \beta X \times \beta X$ via $\overline{(r,s)}$.\qed
\end{lemma}

Using this Lemma, we can conclude that the pair groupoid operations on $\beta X \times \beta X$ restrict to give continuous operations on $G(X)$; we equip $G(X)$ with this induced product and inverse.

Given this construction of $G(X)$ as an extension of the pair groupoid, it is easy to see that the set $X$ is an open \textit{saturated} subset of $\beta X$. In particular, this also means that the Stone-\v{C}ech boundary $\partial\beta X$ is saturated. 
\begin{definition}\label{def:bdry}
The \textit{boundary groupoid} associated to $X$ is the restriction groupoid $G(X)|_{\partial\beta X}$.
\end{definition}

\section{Negative type functions on groupoids}\label{sect:negtype}
The role of positive and conditionally negative type kernels within group theory is well known and plays an important role in studying both anayltic and representation theoretic properties of groups \cite{MR2415834,MR1487204}. These notions were considered for groupoids by Tu in \cite{MR1703305} and we use his conventions here. Let $\G$ be a locally compact, Hausdorff groupoid.

\begin{definition}
A continuous function $F: \G \rightarrow \mathbb{R}$ is said to be of \textit{negative type} if 
\begin{enumerate}
\item $F|_{\G^{(0)}}=0$;
\item $\forall x \in \G, F(x)=F(x^{-1})$;
\item Given $x_{1},...,x_{n} \in \G$ all having the same range and $\sigma_{1},...,\sigma_{n} \in \mathbb{R}$ such that $\sum_{i}\sigma_{i}=0$ we have $\sum_{j,k}\sigma_{j}\sigma_{k}F(x_{j}^{-1}x_{k})\leq 0$.
\end{enumerate}
\end{definition}

The important feature of functions of this type is their connection to a-T-menability for locally compact, $\sigma$-compact groupoids, in fact the following are equivalent by work of Tu \cite{MR1703305}:
\begin{enumerate}
\item There exists a locally proper negative type function on $\G$.
\item There exists a continuous field of Hilbert spaces over $\G^{(0)}$ with a proper affine action of $\G$.
\end{enumerate}
We remark that in the situation we are considering, \'etale groupoids with compact Hausdorff base space, local properness is equivalent to properness. We also remark that an alternative proof of this is given by Renault in Section 2.4 of \cite{Renault2012}.

We can now make precise what we mean by ``$G(X)$ captures the coarse geometry of $X$'' in our particular case:

\begin{theorem}\label{thm:coarseprop}
Let $X$ be a uniformly discrete space with bounded geometry and let $G(X)$ be its coarse groupoid. Then $X$ admits a coarse embedding into Hilbert space if and only if $G(X)$ is a-T-menable \cite[Theorem 5.4]{MR1905840}\qed
\end{theorem}

\section{Main Theorem}
We dedicate this section to proving Theorem \ref{Thm:MT1}. More precisely we construct, given a fibred coarse embedding into Hilbert space, a proper conditionally negative type function on the boundary groupoid $G(X)|_{\partial\beta X}$.

\begin{definition}
Let $\mathcal{F}:=\lbrace K_{R} \rbrace$ be a family of bounded subsets of X. Let $A_{R}$ be the restricted entourage $\Delta_{R}\cap ((X\setminus K_{R}) \times (X\setminus K_{R}))$. 
\end{definition}

\begin{rem}\label{rem:outline}
Let $X$ admit a fibred coarse embedding into Hilbert space.
\begin{enumerate}
\item $A_{R}$ is an entourage with the same corona as $\Delta_{R}$, that is $\overline{A_{R}}\setminus A_{R}=\overline{\Delta_{R}}\setminus \Delta_{R}$. Assume that $X$ admits a fibred coarse embedding into Hilbert space with respect to $\mathcal{F}$. Then for each $R>0$ the following function is defined on $A_{R}$:
\begin{equation*}
k_{R}(x,y)=\Vert t_{x}(x)(s(x))-t_{x}(y)(s(y))\Vert^{2}
\end{equation*}
where $t_{x}:=t_{B_{R}(x)}:H_{x} \rightarrow H$ is the affine isometry provided by the fibred coarse embedding.
\item In the special case that $X$ is a coarse disjoint union constructed from a sequence of finite graphs $\lbrace X_{i} \rbrace$, the $A_{R}$'s defined above can be taken to have the form: $A_{R}= \sqcup_{i\geq i_{R}}\Delta_{R}^{i}$ for some $i\geq i_{R}$ and where $\Delta_{R}^{i}$ denotes the $R-$neighbourhood of the diagonal in $X_{i}$.
\item For $G(X)|_{\partial\beta X} = \bigcup_{R>0} \overline{\Delta_{R}}\setminus \Delta_{R}$ there is a smallest $R>0$ such that $\gamma \in \overline{\Delta_{R}}\setminus \Delta_{R}$. We denote this by $R_{\gamma}$.
\end{enumerate}
\end{rem}

Using this data we wish to construct a proper, negative type function $f$ on $G(X)|_{\partial\beta X}$. The main idea is to extend each $k_{R}$ to a function $\tilde{k}_{R}$ using the universal property of $\beta A_{R} \cong \overline{A_{R}}$. We then show that these extensions $\tilde{k}_{R}$ piece together to form a single function $f$ on $G(X)|_{\partial\beta X}$. The next section gives the technical details of how to do this supposing that the following additional criterion is satisfied.

\begin{definition}\label{def:scaleindep}
A family of kernels $\lbrace k_{R}: A_{R} \rightarrow \mathbb{R}_{+}\rbrace$ is \textit{scale independent} if for every $S>R>0$, we have that $k_{R}|_{A_{R}\cap A_{S}} = k_{S}|_{A_{R}\cap A_{S}}$.   
\end{definition}

This notion will allow us to piece the kernels on every scale together as prescribed above.

\subsection{The technical steps.}
In this section, we give the details of the construction outlined above.

\begin{proposition}\label{prop:wd}
Let $X$ be a uniformly discrete bounded geometry metric space that fibred coarsely embeds into Hilbert space, and suppose additionally that the kernels $k_{R}$ coming from this embedding are scale independent. Then the function $f$ obtained by patching together the $k_{R}$ is well defined and continuous on $G(X)|_{\partial\beta X}$.
\end{proposition}
\begin{proof}
We have two considerations:
\begin{enumerate}
\item Extending $k_{R}$ to $\tilde{k}_{R}$ on the closure of $A_{R}$ for each $R>0$.
\item Each $\tilde{k}_{R}$ pieces together; for any $S>R$ we have that $\tilde{k}_{S}(\gamma)=\tilde{k}_{R}(\gamma)$.
\end{enumerate}
We first prove (1). Under the assumption that $X$ fibred coarsely embeds into Hilbert space, we know that each $k_{R}$ is a bounded function on $A_{R}$, hence extends to a continuous function $\tilde{k}_{R}$ on the Stone-\v{C}ech compactification of $A_{R}$, which in this context is homeomorphic to its closure in $\beta(X \times X)$. 

For the proof of (2) consider the sets $A_{R}$ and $A_{R}\cap A_{S}$ for $S>R$. Using the compatibility properties of a fibred coarse embedding and the scale independence of Definition \ref{def:scaleindep} we can see that the function $k_{R}$, and $k_{S}$ restricted to $A_{R}\cap A_{S}$ agree. Using the observation that for all $R>0$ and $S>R$ the set $A_{R}\cap A_{S}$ has the same Stone-\v{C}ech boundary as $A_{R}$ we can deduce that $\tilde{k}_{R}$ and $\tilde{k}_{S}$ restricted to $\overline{A_{R}}\setminus A_{R}$ agree. Hence we can define, for any $\gamma \in G(X)|_{\partial\beta X}$, $f(\gamma) = \tilde{k}_{R_{\gamma}}(\gamma)$, which is the natural continuous function defined on the union $\cup_{R>0}\overline{A_{R}}$.
\end{proof}

\begin{lemma}\label{lem:MT1-a}
The function $f$ is proper.
\end{lemma}
\begin{proof}
To see $f$ is proper it is enough to prove that the preimage of an interval $[0,r]$ is contained in $\overline{A_{R}}$ for some $R>0$. This is as each interval $[0,r]$ is a closed subset of $\mathbb{R}$, the map $f$ is continuous and hence $f^{-1}([0,r])$ would be a closed subset of a compact set, hence would itself be compact.

We now assume for a contradiction that the preimage of $[0,r]$ contains elements $\gamma$ with the $R_{\gamma}$ defined in Remark \ref{rem:outline}.(4) being arbitrarily large. Then from the definition of $f$ and the fact that $X$ admits a fibred coarse embedding, we can see:
\begin{equation*}
\rho_{-}(R_{\gamma})^{2} \leq f(\gamma) \leq \rho_{+}(R_{\gamma})^{2}.
\end{equation*}
As $\rho_{-}(S)$ tends to infinity as $S$ does, we can find an $S>0$ such that $\rho_{-}(S)>r$. By assumption there exists $\gamma \in f^{-1}([0,r])$ with $R_{\gamma}$ as large as we like: in particular $\rho_{-}(R_{\gamma})>r$, which is impossible. Whence there exists an $R>0$ such that $f^{-1}([0,r]) \subset \overline{A_{R}}$.
\end{proof}

\begin{lemma}\label{lem:MT1-b}
The function $f$ is of negative type.
\end{lemma}
\begin{proof}
This relies on the ideas of \cite[Theorem 5.4]{MR1905840}. Let $\gamma_{1},...,\gamma_{n} \in G(X)|_{\partial\beta X}$ such that $r(\gamma_{1})=...=r(\gamma_{n}):=\omega$ and let $\sigma_{1},...,\sigma_{n} \in \mathbb{R}$ with sum $0$. We need to prove that:
\begin{equation*}
\sum_{i,j}\sigma_{i}\sigma_{j}f(\gamma_{i}^{-1}\gamma_{j}) \leq 0.
\end{equation*}
As there are only finitely many $\gamma_{i}$, there exists a smallest $R>0$ such that each $\gamma_{i}$ and each product $\gamma^{-1}_{j}\gamma_{i}$ are elements of $\overline{\Delta_{R}}$. Let $(x_{i,\lambda},y_{i,\lambda})$ be nets within $\Delta_{R}$ that converge to $\gamma_{i}$ respectively. As the ranges of the $\gamma_{i}$ are all equal without loss of generality we can assume that $y_{\lambda,i}=y_{\lambda}$ is equal in each net. To see this is possible take an approximation of $\omega$ by a net $y_{\lambda}$ and then use the fact that each $\gamma_{i}$ belongs to the closure of the graph of a partial translation $t_{i}$ \cite{MR2007488}. Now use $x_{\lambda,i}=t_{i}^{-1}(y_{\lambda})$: $\gamma_{i} = \lim_{\lambda}(x_{\lambda,i},y_{\lambda})$.

Hence, for each $\lambda$, we know that $(x_{\lambda,i},y_{\lambda})$ and $(x_{\lambda,j},x_{\lambda,i}) \in \Delta_{R}$, and that $x_{\lambda,i}, x_{\lambda,j}\in B_{R}(y_{\lambda})$.

We wish to compute $\sum_{i,j}\sigma_{i}\sigma_{j}k_{R}(x_{\lambda, j},x_{\lambda, i})$ by using an affine isometry to work relative to $y_{\lambda}$: denote by $t_{y_{\lambda}}$ the local trivialisation $t_{B_{R}(y_{\lambda})}$ and by $k_{R}^{y_{\lambda}}$ the kernel defined using the trivialisation $t_{y_{\lambda}}$ as in Remark 17.1. Let $t_{x_{\lambda,j},y_{\lambda}}$ be the unique affine ismometry such that $t_{x_{\lambda,j},y_{\lambda}}t_{y_{\lambda}}=t_{x_{\lambda,j}}$ on the intersections of the corresponding balls, which we know exist from the definition of a fibred coarse embedding.

Hence: 
$$
\sum_{i,j}\sigma_{i}\sigma_{j}k_{R}(x_{\lambda, j},x_{\lambda, i})  =  \sum_{i,j}\sigma_{i}\sigma_{j}\Vert t_{x_{\lambda,j}}(x_{\lambda ,j})(s(x_{\lambda ,j})) - t_{x_{\lambda,j}}(x_{\lambda ,i})(s(x_{\lambda ,i})) \Vert^{2}
$$
$$ =  \sum_{i,j}\sigma_{i}\sigma_{j}\Vert t_{x_{\lambda,j},y_{\lambda}}t_{y_{\lambda}}(x_{\lambda ,j})(s(x_{\lambda ,j})) - t_{x_{\lambda,j},y_{\lambda}}t_{y_{\lambda}}(x_{\lambda ,i})(s(x_{\lambda ,i})) \Vert^{2}\\
$$
$$
 =  \sum_{i,j}\sigma_{i}\sigma_{j}k_{R}^{y_{\lambda}}(x_{\lambda, j},x_{\lambda, i}) \mbox{ as } t_{x_{\lambda,j},y_{\lambda}} \mbox{ is an isometry}.
$$
This reformulation allows us to directly compute using a standard argument (although with worse notation):
\begin{eqnarray*}
&&\sum_{i,j}\sigma_{i}\sigma_{j}k_{R}^{y_{\lambda}}(x_{\lambda, j},x_{\lambda, i}) = \sum_{i,j}\sigma_{i}\sigma_{j}\Vert t_{y_{\lambda}}(x_{\lambda ,j})(s(x_{\lambda ,j})) - t_{y_{\lambda}}(x_{\lambda ,i})(s(x_{\lambda ,i})) \Vert^{2}\\
&&= \sum_{i,j}\sigma_{i}\sigma_{j}(\Vert t_{y_{\lambda}}(x_{\lambda ,j})(s(x_{\lambda ,j})) \Vert^{2} + \Vert t_{y_{\lambda}}(x_{\lambda ,i})(s(x_{\lambda ,i})) \Vert^{2} - 2 \langle t_{y_{\lambda}}(x_{\lambda ,j})(s(x_{\lambda ,j})), t_{y_{\lambda}}(x_{\lambda ,i})(s(x_{\lambda ,i}))\rangle)\\
&&=(\sum_{j}\sigma_{j}\Vert t_{y_{\lambda}}(x_{\lambda ,j})(s(x_{\lambda ,j})) \Vert^{2})(\sum_{i}\sigma_{i})+(\sum_{i}\sigma_{i}\Vert t_{y_{\lambda}}(x_{\lambda ,i})(s(x_{\lambda ,i})) \Vert^{2})(\sum_{j}\sigma_{j})\\
&&  -2\langle \sum_{j}\sigma_{j}t_{y_{\lambda}}(x_{\lambda ,j})(s(x_{\lambda ,j})),\sum_{i}\sigma_{i}t_{y_{\lambda}}(x_{\lambda ,i})(s(x_{\lambda ,i}))\rangle \leq 0.
\end{eqnarray*}
This holds for each $\lambda$ in the net. Taking a limit in $\lambda$:
\begin{equation*}
\sum_{i,j}\sigma_{i}\sigma_{j}f(\gamma_{i}^{-1}\gamma_{j})=\sum_{i,j}\sigma_{j}\sigma_{i}\lim_{\lambda}k_{R}^{y_{\lambda}}(x_{\lambda ,j},x_{\lambda, i}) \leq 0
\end{equation*}
\end{proof}

The preceding Lemmas will now be used to prove the following Theorem.

\begin{theorem}\label{thm:MT1-a}
Let $X=\sqcup_{i}X_{i}$ be a coarse disjoint union of finite metric spaces such that the space $X$ is uniformly discrete and has bounded geometry. If $X$ admits a fibred coarse embedding into Hilbert space then the boundary groupoid $G(X)|_{\partial\beta X}$ is a-T-menable. 
\end{theorem}
\begin{proof}
By assumption the fibred coarse embedding can be chosen such that the family of kernels satisfies the restriction property of Definition \ref{def:scaleindep} (this follows from the fact that such spaces admit the nicer Definition \ref{def:FCE2}, which is Definition 5.1 of \cite{MR3116568}). The Lemmas above now combine to prove the result. 
\end{proof}

The remainder of this section is to show how Theorem \ref{thm:MT1-a} can be adapted to the situation that $X$ is is not a coarse disjoint union but is fibred coarsely embeddable into Hilbert space.

\subsection{What happens if the space is not a disjoint union of finite metric spaces?}\label{sect:whathappens}

In general, the property of scale independence of Definition \ref{def:scaleindep} will not hold for an arbitrary fibred coarse embedding (at least not obviously). However, in this section we outline a method of ``coarse'' decomposition that will still allow us to prove the analogue of Theorem \ref{thm:MT1-a} for general uniformly discrete spaces with bounded geometry. We will use \cite[Chapter 3]{MR2562146} as a common reference for results concerning the relationship between coarse embeddings and kernels in this section.

\begin{definition}
Let $Z$ be a uniformly discrete bounded geometry metric space and let $z_{0} \in Z$. Then an \textit{annular decomposition} of $Z$ is a covering $Y_{0}\cup Y_{1}=Z$, where $Y_{0}, Y_{1}$ and $Y_{0}\cap Y_{1}$ are unions of annular regions around $z_{0}$. 
\end{definition}
The next Lemma, the details of which appear in \cite{MR3116568}, illustrates it is always possible to find a decomposition that is very well spaced, i.e a decomposition in which each of $Y_{0},Y_{1}$ and $Y_{0}\cap Y_{1}$ is a coarse disjoint union of finite metric spaces. 

\begin{lemma}\label{lem:annular}
Let $Z$ be a uniformly discrete bounded geometry metric space. Then $Z$ admits an annular decomposition that is coarsely excisive.
\end{lemma}
\begin{proof}
We give the construction here and the remainder is implicit from the proof of Theorem 1.1, from page 19 of \cite{MR3116568}. Let $z_{0} \in Z$, then for $n\in \mathbb{N}$ let \begin{equation*}
Z_{n}=\lbrace z \in X| n^{3}-n \leq d(z,z_{0}) \leq (n+1)^{3}+(n+1)\rbrace.
\end{equation*}
Now, set $Y_{0}$ to be the union over even $n$ and $Y_{1}$ to be the union over odd $n$. 
\end{proof}
So, for any uniformly discrete bounded geometry metric space $Z$ we have two decompositions:
\begin{enumerate}
\item $Z=Y_{0}\cup Y_{1}$, where the two pieces are themselves unions of finite (annular) metric spaces;
\item $Z=\cup_{n} Z_{n}$ , where each $Z_{n}$ are finite but not of uniformly bounded cardinality in $n$. 
\end{enumerate}
 
We can still use the $Z_{n}$ as a covering of $X$ to define \textit{partition of unity} and because they form the sets $Y_{0}$ and $Y_{1}$, we can use the favourable notion of fibred coarse embedding from Definition \ref{def:FCE2} when we work with them. To get a partition of unity we use Proposition 4.1 of \cite{MR2364071}:

\begin{proposition}\label{prop:pou}
Let $\mathcal{U}$ be a cover of a metric space $X$ with multiplicity at most $k+1$ and Lebesgue number $L>0$. For $U \in \mathcal{U}$ define:
\begin{equation*}
\phi_{U}(x)= \frac{d(x,X \setminus U)}{\sum_{V\in \mathcal{U}}d(x,X \setminus V)}
\end{equation*}
then $\lbrace \phi_{U} \rbrace_{U\in \mathcal{U}}$ is a partition of unity on $X$ subordinated to the cover $\mathcal{U}$. Moreover, each $\phi_{U}$ satisfies:
\begin{equation*}
\vert \phi_{U}(x)- \phi_{U}(y)\vert \leq \frac{2k+3}{L}d(x,y) \mbox{ for every } x,y \in X.
\end{equation*}\qed
\end{proposition}

Let $X$ be a uniformly discrete space, which we will later take to be fibred coarsely embeddable. We will now attempt to apply Proposition \ref{prop:pou} to the covering of $X$ obtained by applying Lemma \ref{lem:annular}, that is $X=\cup_{n} Z_{n}$. Define, $X_{> n}$ to be the union $\cup_{m > n} Z_{m} \subset X$.

We observe that by the construction in Lemma \ref{lem:annular} the multiplicity of the covering by the $Z_{n}$ is $2$. Furthermore, we can control the Lebesgue number of this cover away from a bounded subset using the following observation:

\begin{lemma}\label{lem:estimate}
Let $x,y$ belong to the symmetric difference of $Z_{n}$ and $Z_{n+1}$. If $x\in Z_{n}$ and $y\in Z_{n+1}$ then $d(x,y)\geq 2(n+1)$.
\end{lemma} 
\begin{proof}
If $x\in Z_{n}$ and $y\in Z_{n+1}$ then by the definition of $Z_{n}$ and the fact that $x\not \in Z_{n+1}$ tells us that $d(z_{0},x)\leq (n+1)^{3}-(n+1)$. Similiarly, as $y\in Z_{n+1}$ but $y\in Z_{n}$, we know that $d(z_{0},y)\geq (n+1)^{3}+(n+1)$. Combining these inequalities with the triangle inequality gives the desired estimate.
\end{proof}

An immediate consequence of this is that for $L>0$ there is an $n_{L}$ such that the covering restricted to $X_{>n_{L}}$ has Lebesgue number $\leq L$. To see this, consider the situation that there is a point $x\in Z_{n}\subset Z$ such that the ball $B_{L}(x)$ is not contained in some element $Z_{m}$ of the covering. As the multiplicity of the cover is $2$, it is immediate that the ball cannot be contained in any sets other than $Z_{n-1}$, $Z_{n}$ or $Z_{n+1}$. We can conclude that, without loss of generality, there is a $y\in Z_{n+1}$ such that there is a point $y\in Z_{n+1}$ that does not belong to $Z_{n}$ but does belong to the ball $B_{L}(x)$. If $x \not \in Z_{n+1}$, then Lemma \ref{lem:estimate} implies that $L \geq 2(n+1)$. So by choosing $n_{L}$ to be the least integer above $\frac{L}{2}-1$ we see that this situation cannot occur. What remains is the case $x\in Z_{n+1}$, but a mirror of this argument reduces this to the previous case. Hence on the space $X_{> n_{L}}$ the covering has Lesbegue number at least $L$.

Suppose now that $X$ fibred coarsely embeds into Hilbert space with respect to a family $\lbrace K_{R} \rbrace$. Then for every $R>0$ there is a smallest natural number such that the bounded set $K_{R}$ is contained within some finite union of annuli. Denote this by $m_{R}$. Observe also that the above statement can be modified to take into account the two bounded sets that will occur after we break $X$ into $Y_{0}$ and $Y_{1}$. In this case, we define $m_{R}$ to be the maximum of the natural numbers required on each piece.

This allows us to control, up to scale, the asymptotic behaviour of kernel functions on the spaces $Y_{0}$ and $Y_{1}$ when gluing them back together. To make this formal, we need first to recall the result of Schoenberg \cite{MR1501980}, where the statement here is taken from \cite{MR2562146}:

\begin{theorem}\label{thm:schoenberg}
Let $X$ be a set and let $k$ be a negative type kernel on $X$. Then for every $t>0$, the kernel $F_{t}(x,y):=e^{-tk(x,y)}$ is of positive type.\qed
\end{theorem}

On each of the spaces $Y_{i}$ we have a family of scale independent negative type kernels $\lbrace k_{i,R}\rbrace_{R}$ defined on $A_{R}$ that is well controlled by a common pair of functions $\rho_{\pm}$. By applying Theorem \ref{thm:schoenberg}, we now have a family, controlled by $R$ and $t$, of scale independent kernels of positive type that are also well controlled (c.f \cite{MR2562146} Chapter 3 for a full treatment). Recall that a kernel $k$ is said to have $(R,\epsilon)$-variation if whenever $(x,y)\in \Delta_{R}$, we have that $\vert 1 - k(x,y)\vert \leq \epsilon$.

The following result collects these ideas and is a modification of an argument shared with the author by Rufus Willett. 
\begin{lemma}\label{lem:complicatedsch}
Suppose $X$ is fibred coarsely embeddable into Hilbert space. Then there is a family of normalised positive type kernels $\lbrace F_{R,t}:A_{R}\rightarrow \mathbb{R}\rbrace_{R,t}$ such that the following are satisfied:
\begin{itemize}
\item $\lbrace F_{R} \rbrace_{R}$ is a scale independent family for every $t$;
\item for every $R,\epsilon>0$ there is a bounded subset $D$ that contains $K_{R}$ and a real number $t$ such that there is a positive type kernel $F_{R,t}$ with $(R,\epsilon)$-variation when restricted to $A_{R}\cap (X\setminus D)^{2}$;
\item for all $t,\epsilon > 0$ there exists $S>0$ such that for all $R>S$ and all $(x,y)\in A_{R}$ that do not belong to $A_{S}$ we have $\vert F_{R,t}(x,y)\vert < \epsilon$.
\end{itemize}
\end{lemma}
\begin{proof}
The proof proceeds in a few steps:
\begin{enumerate}
\item We define, from Proposition \ref{prop:pou}, an $\ell^{2}$-partition of unity. 
\item Using this partition of unity for the cover $\mathcal{U}=\lbrace Z_{n} \rbrace$ and the families of kernels on $Y_{i}$ provided by Theorem \ref{thm:schoenberg} and the assumption of fibred coarse embeddability we define a family of new kernels for $X$. This family will be scale independent.
\item We show that this kernel can be constructed such that for any $R$ and $\epsilon$, there is a natural number $o_{R,\epsilon}$ and a $t$ such that the kernel $F_{R,t}$ restricted to the set $X_{>o_{R,\epsilon}}$ has $(R,\epsilon)$-variation.
\item Finally, we prove the asymptotic control estimates using observations concerning the control functions that come from the original fibred coarse embedding.
\end{enumerate}
Let $n_{R}$ and $m_{R}$ be defined as during the discussion prior to the statement of the Lemma. We now work through the steps outlined above:

\begin{enumerate}
\item Consider the new partition of unity defined using $\Phi_{n}=\sqrt{\phi_{X_{n}}}$. This will have $\ell^{2}$-sum equal to $1$ for each $x$, and satisfies, for any pair $x,y\in X$:
\begin{equation*}
\vert \Phi_{n}(x)  - \Phi_{n}(y) \vert \leq \sqrt{\frac{5}{L}d(x,y)}.
\end{equation*}
\item Now, on each scale we restrict this partition of unity to the subspace $X\setminus D_{R}$ and observe that for all but finitely many $x \in X$, the sum: $\sum_{n>\max(n_{R},m_{R})} \Phi_{n}(x) = \sum_{n} \Phi_{n}(x)$. Let $F_{R,t}^{[n]}$ be the kernel output by Schoenbergs theorem that arises from the kernel $k_{[n],R}$, where $[n]$ represents the class of $n$ modulo $2$, and consider for each $R>0$ the kernel $F_{R,t}$ defined pointwise by:
\begin{equation}\label{eqn:kernel}
F_{R,t}(x,y) = \sum_{n}\Phi_{n}(x)F_{R,t}^{[n]}(x,y)\Phi_{n}(y)
\end{equation}

By construction this family of kernels is scale independent for each $t$ as each $F_{R,t}^{i}$ are, and the same partition of unity is used and it was constructed globally. This completes Step 2.

\item By combining observations from above and estimates from Theorem 3.2.8(2) in \cite{MR2562146} we can see that for every $R,\epsilon>0$ there exists an $o_{R,\epsilon}$ such that the covering restricted to $X_{>o_{R,\epsilon}}$ has Lebesgue number $S:=\frac{180R}{\epsilon^{2}}$ and $K_{R} \subseteq \cup_{n\leq o_{R,\epsilon}}Z_{n}:=D$.

Let $t=\frac{\epsilon}{3(1+\rho_{+}(R)^{2})}$. Then each of the kernels $F_{R,t}^{i}$ satisfies 
\begin{equation*}
1-F^{i}_{R,t}(x,y) \leq 1- e^{-t\rho_{+}(d(x,y))}\leq 1-e^{\frac{-\epsilon}{3}}\leq \frac{\epsilon}{3}
\end{equation*}
for every $x,y \in A_{R}$. Now, consider the following for every $x,y\in A_{R}\cap (X\setminus D)^{2}$:
\begin{eqnarray*}
1-F_{R,t}(x,y) & = & \sum_{n}\Phi_{n}(x)^{2} - \Phi_{n}(x)F_{R,t}^{[n]}(x,y)\Phi_{n}(y)\\
 & = & \Phi_{n}(x)^{2} + \Phi_{n+1}(x)^{2} - \Phi_{n}(x)F_{R,t}^{[n]}(x,y)\Phi_{n}(y) - \Phi_{n+1}(x)F_{R,t}^{[n+1]}(x,y)\Phi_{n+1}(y)\\
 & \leq & \Phi_{n}(x)\Phi_{n}(y) + \Phi_{n}(x)(\Phi_{n}(x) - \Phi_{n}(y))\\
 & & + \Phi_{n+1}(x)\Phi_{n+1}(y) + \Phi_{n+1}(x)(\Phi_{n+1}(x) -\Phi_{n+1}(y)) \\
 & & - \Phi_{n}(x)F_{R,t}^{[n]}(x,y)\Phi_{n}(y) - \Phi_{n+1}(x)F_{R,t}^{[n+1]}(x,y)\Phi_{n+1}(y) \\
  & < & \Phi_{n}(x)\Phi_{n}(y)(1-F_{R,t}^{[n]}(x,y)) + \Phi_{n+1}(x)\Phi_{n+1}(y)(1-F_{R,t}^{[n+1]}(x,y)) + \frac{\epsilon}{3}.
\end{eqnarray*}
Where $[n]$ denotes $n$ modulo $2$, $n$ depends only on $(x,y)$ and the $\frac{\epsilon}{3}$ term comes from the choice of $o_{R,\epsilon}$ such that each $\vert \Phi_{n}(x) - \Phi_{n}(y) \vert < \sqrt{\frac{5R}{S}} = \frac{\epsilon}{6}$. To complete the proof, observe also that $\vert \Phi_{n}(x)\Phi_{n}(y)(1-F_{R,t}^{[n]}(x,y)) \vert \leq \frac{\epsilon}{3}$ by the choice of $S$ and $t$.

\item We now use the fact that these kernels are scale independent; the $(R,\epsilon)$-variation property above holds only on a large piece of $A_{R}$ and outside this piece the kernels $F^{i}_{R,t}$ satisfy $F^{i}_{R,t}(x,y)= e^{-tk_{R}(x,y)} \leq e^{-t\rho_{-}(d(x,y))}$, which clearly tends to $0$ as $d(x,y)\rightarrow \infty$. Patching this together, we see that $F_{R,t}$ will also tend to $0$ as $d(x,y)\rightarrow \infty$, and so for fixed $t$ and $\epsilon$ we pick $S$ large enough such that $e^{-t\rho_{-}}(S)$ is less than $\epsilon$, which proves point 3 in the claim.
\end{enumerate}
\end{proof}

The important point of Lemma \ref{lem:complicatedsch} is that the family of kernels $\lbrace F_{R,t}\rbrace_{R}$ is extendable to the boundary groupoid of $X$ for every $t$, so we can now prove an analogue of Theorem 3.2.8(4) of \cite{MR2562146} at infinity using these families.

\begin{theorem}
Let $X$ be a uniformly discrete metric space that admits a fibred coarse embedding into Hilbert space. Then the boundary groupoid $G(X)|_{\partial\beta X}$ is a-T-menable.
\end{theorem}
\begin{proof}
We proceed as in Theorem 3.2.8 of \cite{MR2562146}.  
\begin{enumerate}
\item For every $n\in \mathbb{N}$ define, using Lemma \ref{lem:complicatedsch}, a family of kernels $F_{n}$ that satisfy $1-F_{n}(x,y) < 2^{-n}$ whenever $d(x,y)\leq n$. 
\item We now consider the extensions of $F_{n}$ to the boundary groupoid $G(X)|_{\partial\beta X}$, which is possible as they are each scale independent.
\item Define $k(\gamma)=\sum_{n=1}^{\infty}(1-\tilde{F_{n}}(\gamma))$. This is the limit of the partial sums defined on the boundary of the closure of $(X\setminus D_{n,2^{-n}})^{2}$. This converges as for any $\gamma \in G(X)|_{\partial\beta X}$ as there are only finitely many $n$ for which $\gamma \not \in \overline{\Delta_{n}}\setminus\Delta_{n}$ and let $n_{\gamma}$ denote the largest such $n$. Now the tail of this sum will be bounded above by $\sum_{n>n_{\gamma}}^{\infty}2^{-n}$. The function $k$ has negative type, as each term in the sum does and this property is closed under pointwise limits and positive sums (The extension point follows from nothing other than the argument for Lemma \ref{lem:MT1-b}).
\item It is now enough to show that this is proper (and we do this by showing it is controlled by two functions that tend to infinity $\tau_{\pm}$). Observe that for each $\gamma \in \overline{\Delta_{R}}\setminus\Delta_{R}$ we have:
\begin{equation*}
k(\gamma) \leq \sum_{n=1}^{\lfloor R \rfloor}1 + \sum_{n=\lfloor R +1\rfloor}^{\infty}2^{-n} \leq R+1
\end{equation*}
so we can take $\tau_{+}(R)=R+1$. For the lower bound, let $h_{n}(R):=e^{-t_{n}\rho_{-}(R)^{2}}$, where $t_{n}$ is the parameter that realises the kernel $F_{n}$. Then for each $N \in \mathbb{N}$, there exists an $S_{N}$ such that for all $n \leq N$ and $S>S_{N}$, we have $h_{n}(S)<\frac{1}{2}$. We can find such an $S_{N}$ as each $h_{n}$ tends to $0$ at infinity. Now, for every $\gamma \in G(X)|_{\partial\beta X}$ that is a limit of pairs in $A^{c}_{S_{N}}$, we have:
\begin{equation*}
k(\gamma) \geq \sum_{n=1}^{N}1-F_{n}(\gamma) \geq \sum_{n=1}^{N}1-h_{n}(S_{N}) \geq \frac{N}{2}.
\end{equation*}
Finally, choose:
\begin{equation*}
\tau_{-}(s) = \frac{1}{2}\max\lbrace N | s \geq s_{N} \rbrace.
\end{equation*}
With this choice we have $\tau_{-}(R) \leq k(\overline{\Delta_{R}}\setminus\Delta_{R}) \leq \tau_{+}(R)$. Hence, $k$ is a proper, negative type function on $G(X)|_{\partial\beta X}$.
\end{enumerate}
\end{proof}

\section{Applications}\label{sect:apps}

\subsection{The boundary coarse Baum-Connes conjecture and the coarse Novikov conjecture}\label{Sect:apps}
Throughout this section let $A_{\partial}$ denote the quotient $C^{*}$-algebra $l^{\infty}(X,\mathcal{K})/C_{0}(X,\mathcal{K})$. Recall from \cite{mypub1} the boundary coarse Baum-Connes conjecture.
\begin{conjecture1} [Boundary Coarse Baum-Connes Conjecture]
Let $X$ be a uniformly discrete bounded geometry metric space. Then the assembly map:
\begin{equation*}
\mu_{bdry}:K_{*}^{top}(G(X)|_{\partial\beta X}, A_{\partial}) \rightarrow K_{*}(A_{\partial}\rtimes_{r}G(X)|_{\partial\beta X})
\end{equation*}
is an isomorphism.
\end{conjecture1}

This conjecture also has a maximal form \cite[Section 4]{mypub1} that is equivalent to the maximal coarse Baum-Connes conjecture at infinity from \cite{MR3116568}. If $X$ is a uniformly discrete bounded geometry metric space we can see that the algebra at infinity defined in \cite{MR3116568} and the groupoid crossed product algebra $A_{\partial}\rtimes_{m}G(X)|_{\partial\beta X}$ are isomorphic. Proceeding via this conjecture we can appeal to the machinery of Tu \cite{MR1703305} concerning $\sigma$-compact, locally compact a-T-menable groupoids to conclude results about the (maximal) coarse Baum-Connes conjecture. 

In particular, we can use this conjecture and homological algebra to conclude Theorem \ref{Thm:FCEMR}.

\begin{theorem}\cite[Theorem 1.1]{MR3116568}
Let $X$ be a uniformly discrete space with bounded geometry that fibred coarse embeds into Hilbert space. Then the maximal coarse Baum-Connes assembly map is an isomorphism for $X$.
\end{theorem}
\begin{proof}
We have a short exact sequence of maximal groupoid $C^{*}$-algebras:
\begin{equation*}
0 \rightarrow \mathcal{K} \rightarrow C^{*}_{m}(G(X)) \rightarrow C^{*}_{m}(G(X)|_{\partial\beta X}) \rightarrow 0.
\end{equation*}
This gives us the following diagram, arising from the long exact sequence in K-theory and suitable Baum-Connes conjectures (omitting the coefficients):
\begin{equation*}
\xymatrix@=0.7em{
K_{1}(C^{*}(G(X)|_{\partial\beta X})) \ar[r] & K_{0}(\mathcal{K}) \ar[r]& K_{0}(C^{*}_{m}(G(X))) \ar[r]& K_{0}(C^{*}_{m}(G(X)|_{\partial\beta X}))\ar[r] & K_{1}(\mathcal{K})  \\
K_{1}^{top}(G(X)|_{\partial\beta X}) \ar[r] \ar[u]& K_{0}^{top}(X \times X) \ar[r]\ar[u]^{\ucong}& K_{0}^{top}(G(X)) \ar[r]\ar[u]& K_{0}^{top}(G(X)|_{\partial\beta X}) \ar[r]\ar[u]^{\mu_{bdry}}& K_{1}^{top}(X \times X)\ar[u]^{\ucong}
}
\end{equation*}
By Corollary \ref{Thm:MT1} the maximal boundary assembly map is an isomorphism. The result now follows from the Five lemma.
\end{proof}

To understand the reduced assembly map requires a more delicate approach.

\begin{definition}
Let $X$ be as above. We say that $X$ has an \textit{infinite coarse component} if there exists $E \in \mathcal{E}$ such that $P_{E}(X)$, the Rips complex over $E$, has an unbounded connected component. Otherwise we say that $X$ \textit{only has finite coarse components}. In the metric coarse structure, this condition becomes: if there exists an $R>0$ such that $P_{R}(X)$ has an unbounded connected component.
\end{definition}

\begin{example}(Space of Graphs)
Let $\lbrace X_{i} \rbrace_{i \in \mathbb{N}}$ be a sequence of finite graphs such that $\vert X_{i} \vert \rightarrow \infty$ in $i$. Then we can form the \textit{coarse disjoint union} (or \textit{space of graphs}) $X$ with underlying set $\sqcup X_{i}$, as in Definition \ref{def:coarselydisconnected}. The space $X$ is then a prototypical example of a space with only finite coarse components.
\end{example}

\begin{lemma}\label{lem:zandi}
Let $X$ be a uniformly discrete bounded geometry metric space, $i: \mathcal{K} \hookrightarrow C^{*}X$ be the inclusion of the compact operators into the Roe algebra of $X$ and $j:X\times X \rightarrow G(X)$ be the canonical inclusion. Then:
\begin{enumerate}
\item If $X$ has an infinite coarse component then the map from $K^{top}(X \times X) \rightarrow K^{top}(G(X))$ is the $0$ map.
\item If $X$ is a space with only finite coarse components then $i_{*}:K_{*}(\mathcal{K}) \rightarrow K_{*}(C^{*}X)$ is injective.
\end{enumerate}
\end{lemma}
\begin{proof}
Claim (1) relies on an argument involving groupoid equivariant KK-theory \cite{MR1656031,MR1686846,MR1798599}. We remark that it is useful in this instance to decompose $G(X)=\beta X \rtimes \G_{A}$, where $\G_{A}$ is some second countable, Hausdorff \'etale groupoid that contains the pair groupoid $X \times X$ as an open subgroupoid. This is possible by Lemma 3.3b) of \cite{MR1905840}. We then identify $KK_{\G_{A}}$ with $KK_{G(X)}$ where appropriate using Theorem 3.8 from \cite{cbcag2}. 

Let $Z=P_{E}(\G_{A})$, a proper $\G_{A}$-space, and $U$ be the preimage of $X$ under the anchoring map and $Z_{x}$ be the preimage of the singleton $x \in X$. We consider the following diagram, with notation following \cite{MR1905840,MR1656031}:
\begin{equation*}
\xymatrix{KK_{\G_{A}}(C_{0}(U),C_{0}(X,\mathcal{K})) \ar[r]_{1} & KK_{\G_{A}}(C_{0}(Z),\ell^{\infty}(X,\mathcal{K})) \ar[r]^{\cong}_{4}  & KK(C_{0}(Z_{x}),\mathcal{K}) \\
& KK(\mathbb{C},\mathbb{C}) \ar_{2}[lu]\ar[u]_{3}\ar[ur]_{5}&
}
\end{equation*}
with maps at the level of cycles given by:
\begin{eqnarray*}
1 & : & (E,\phi , F) \mapsto (E\otimes_{\alpha}\ell^{\infty}(X,\mathcal{K}),\phi \otimes 1 , F \otimes 1)\\
2 & : & (\mathbb{C},1,0) \mapsto (C_{0}(X,\mathcal{K}),\tilde{\psi},0)\\
3 & : & (\mathbb{C},1,0) \mapsto (\ell^{\infty}(X,\mathcal{K}),\tilde{\psi},0)\\
4 & : & (E,\phi , F) \mapsto (E|_{x},\phi|_{x},F|_{x})\\
5 & : & (\mathbb{C},1,0) \mapsto (\mathcal{K},\psi,0)
\end{eqnarray*}
Some remarks now about these maps:
\begin{itemize}
\item Map $1$ is the map induced by extending a cycle from $A \triangleleft B$ to $B$ using the interior tensor product construction given on page 38 of \cite{MR1325694}. We apply this construction to the natural inclusion map, denoted $\alpha: C_{0}(X,\mathcal{K}) \subseteq \ell^{\infty}(X,\mathcal{K})$; where $\ell^{\infty}(X,\mathcal{K})$ is considered as compact adjointable operators over itself considered as a Hilbert $C^{*}$-module. To construct the left action of $C_{0}(Z)$ on $E\otimes_{\alpha}\ell^{\infty}(X,\mathcal{K})$ we take the composition of $\psi: C_{0}(U) \rightarrow \mathbb{B}(E)$ and $\alpha_{*}:\mathbb{B}(\mathcal{E}) \rightarrow \mathbb{B}(E\otimes_{\alpha}\ell^{\infty}(X,\mathcal{K}))$, then extend them using Proposition 2.1 from \cite{MR1325694}. This has the desired properties as $\alpha_{*}$ maps the compacts on $E$ into the compacts on the interior tensor product over $\alpha$ in this instance (by Lemma 4.6 in \cite{MR1325694}).
\item $4$ is the restriction map induced in equivariant KK-theory by the inclusion of a point $\lbrace x \rbrace$ as an open subgroupoid of $G(X)$. This is an isomorphism by Lemma 4.7 from \cite{MR1905840}.
\item  We recall the construction of the natural inverse to $4$ defined in Lemma 4.7 of \cite{MR1905840} denoted here by $\Xi$. Given a cycle $(E,\phi,F)$, we can construct a cycle for $KK_{\G_{A}}(C_{0}(Z),\ell^{\infty}(X,\mathcal{K}))$ by taking the module $\ell^{\infty}(X,\mathcal{K})$, and the representation and operator, $\tilde{\phi}$ and $\tilde{F}$, as follows:
\begin{equation*}
(\tilde{\phi}(a)\xi)(x):=\phi(a(x))\xi(x), \mbox{ and } (\tilde{F}\xi)(x):=F\xi(x).
\end{equation*}

We now define $\Xi$ to be the map that sends the cycle $(E,\phi,F)$ to $(\ell^{\infty}(X,\mathcal{K}), \tilde{\phi}, \tilde{F})$. This map is the same as the map defined in the proof of Lemma 4.7 of \cite{MR1905840} after stabilising the cycles by the degenerate standard cycle $(\widehat{H},0,0)$ (where $\widehat{H}$ is the standard Hilbert module over $\ell^{\infty}(X,\mathcal{K})$).

\item We remark that $5$ is the map induced on K-homology by including a point into $X$. In particular it maps the generating cycle $(\mathbb{C},1,0)$ to $(\mathcal{K},\psi,0)$, where $\psi: C_{0}(Z) \rightarrow \mathcal{K}$ constructed as follows: let $q$ be a rank one projection in $\mathcal{K}$ and fix a point $x \in X$, then define $\psi(f)= f(x)q$. This turns the triple $(\mathcal{K},\psi,0)$ into a Kasparov cycle.

\item Combining these last two points it is possible to describe the maps $2$,$3$ and their compatibility with $5$: We observe that $3$ maps $(\mathbb{C},1,0)$ to $\Xi\circ 5$, which is $(\ell^{\infty}(X,\mathcal{K}),\tilde{\psi},0)$. We can now take $2$ similarly as $\tilde{\psi}$ restricts to $C_{0}(X,\mathcal{K})$ and so we can define $2$ by mapping $(\mathbb{C},1,0)$ to $(C_{0}(X,\mathcal{K}),\tilde{\psi},0)$. This cycle maps using $1$ to $(\ell^{\infty}(X,\mathcal{K}),\tilde{\psi},0)$ (this is seen by considering the interior tensor product defined using an inclusion).
\item We remark also that by the points above it is clear that this diagram commutes.
\end{itemize}

We now take the limit through the directed set of entourages $E \in \mathcal{E}$:
\begin{equation*}
\xymatrix{K^{top}(X \times X, C_{0}(X,\mathcal{K})) \ar[r] & K^{top}(\G_{A},\ell^{\infty}(X,\mathcal{K})) \ar[r]^{\cong}  & KX_{*}(X) \\
& KK(\mathbb{C},\mathbb{C}) \ar^{\cong}[lu]\ar[u]\ar[ur]&
}
\end{equation*}
where Theorem 3.8 \cite{cbcag2} allows us to identify $K^{top}_{*}(X\times X , C_{0}(X,\mathcal{K}))$ with $K^{top}_{*}(\G_{A}, C_{0}(X,\mathcal{K}))$ and $K_{*}^{top}(G(X), \ell^{\infty}(X,\mathcal{K}))$ with $K^{top}_{*}(\G_{A},\ell^{\infty}(X,\mathcal{K}))$.

After this identification, we can observe that map $2$ will induce an isomorphism in the limit at the level of KK-groups (as identified in the diagram). To verify this we first remark that $X$, as an $(X\times X)$-space, is an example of a cocompact classifying space for proper actions and so it is enough to compute the group $KK_{X\times X}(C_{0}(X),C_{0}(X,\mathcal{K}))$.

We do this using the restriction map induced by including $\lbrace x \rbrace$ into $X \times X$, which splits the map $2$ at the level of cycles. Then a straight line homotopy of operators shows that every cycle for $KK_{X\times X}(C_{0}(X),C_{0}(X,\mathcal{K}))$ is homotopic to a cycle that is constant in each fiber. This shows $2$ is an isomorphism in the limit.

Finally, armed with this diagram we can then conclude the result by observing that having an infinite coarse component in $X$ implies that including a point via $5$ factors through the coarse K-homology of a ray and hence is the zero map.

We now prove Claim (2). Let $X$ have only finite coarse components and let $\mathcal{E}$ denote the metric coarse structure on $X$. Now for every entourage $E \in \mathcal{E}$, we can decompose $E$ as $\sqcup_{i \in I} E_{i}$, where each $E_{i}$ is finite, in the following way. Consider the Rips complex $P_{E}X$. We know this decomposes into countably many finite, disjoint complexes $P_{E,i}X$ as we have assumed that we have no ray in any such Rips complex $P_{E}X$. By considering the $0$-skeleton of this complex, which we denote by $X_{i}$, we obtain a disjoint decomposition $X= \sqcup_{i \in I} X_{i}$, and a further decomposition of $E=\sqcup E_{i}$, where $E_{i}:=E \cap (X_{i}\times X_{i})$.

We now consider the coarse structure generated by $E$, which is a substructure of $\mathcal{E}$. Denote the algebraic Roe algebra $\mathbb{C}_{E}[X]$ of finite propagation operators in the coarse structure generated by $E$ and by $C^{*}_{E}X$ its closure as bounded operators on $\ell^{2}(X)$. Now we observe $C^{*}_{E}X$ is the closure of elements supported precisely in each finite $E_{i}$, that is the uniform part of the product $\prod_{i}C^{*}_{E_{i}}X$. Hence it sits naturally in the following diagram:
\begin{equation*}
0 \rightarrow \bigoplus_{i} M_{n_{i}} \rightarrow C^{*}_{E}X_{\infty} \rightarrow \prod_{i}M_{n_{i}}
\end{equation*}
where $n_{i}$ is the size of the $E_{i}$ component of the zero skeleton of the Rips complex $P_{E}X$. We remark also that the sum $\bigoplus_{i} M_{n_{i}}$  is the intersection of the compact operators $\mathcal{K}(\ell^{2}(X))$ with $C^{*}_{E}X$. Denote this ideal, to keep track of objects clearly in a limiting argument, by $\mathcal{K}_{E}$.

The inclusion of the sum into the product induces an injection at the level of K-theory; so, for each $E\in \mathcal{E}$, we now have an injection on K-theory induced by the map $\mathcal{K}_{E}$ into $C^{*}_{E}X$. To complete the proof consider the directed system of diagrams $\lbrace \mathcal{K}_{E}\rightarrow C^{*}_{E}X\rbrace_{E \in \mathcal{E}}$, which converge to the diagram $\mathcal{K} \rightarrow C^{*}X$. These induce natural K-theory maps, which are injections for every $E\in \mathcal{E}$. Whence the map induced by the inclusion $\mathcal{K} \rightarrow C^{*}X$ is injective at the level of K-theory.
\end{proof}

We now have the tools in Lemma \ref{lem:zandi} to prove the following Theorem:

\begin{theorem}\label{thm:mcor1}
Let $X$ be a uniformly discrete bounded geometry metric space that admits a fibred coarse embedding into Hilbert space. Then the coarse Baum-Connes assembly map for $X$ is injective.
\end{theorem}
\begin{proof}
For a uniformly discrete space $X$ with bounded geometry we know that the ghost ideal $I_{G}$ fits into the sequence:
\begin{equation*}
0 \rightarrow I_{G} \rightarrow C^{*}X \rightarrow A_{\partial}\rtimes G(X)|_{\partial\beta X} \rightarrow 0.
\end{equation*}
This addition makes the sequence exact and thus it gives rise to the ladder in K-theory and K-homology where the rungs are the assembly maps defined in \cite{mypub1} (omitting coefficients):
$$
\xymatrix@=1em{
\ar[r] & K_{1}(C^{*}_{r}(G(X)|_{\partial\beta X}) \ar[r] & K_{0}(I_{G}) \ar[r]& K_{0}(C^{*}_{r}(G(X))) \ar[r]& K_{0}(C^{*}_{r}(G(X)|_{\partial\beta X})\ar[r] & K_{1}(I_{G}) \ar[r] & \\
\ar[r] & K_{1}^{top}(G(X)|_{\partial\beta X}) \ar[r] \ar[u]^{\mu^{bdry}}& K^{top}(X \times X) \ar[r]_{1}\ar[u]_{2}& K_{0}^{top}(G(X)) \ar[r]\ar[u]^{\mu}& K_{0}^{top}(G(X)|_{\partial\beta X}) \ar[r]\ar[u]^{\mu^{bdry}}& 0\ar[u]^{0} \ar[r] &
}
$$
Corollary \ref{Thm:MT1} allows us to conclude that $\mu^{bdry}$ is an isomorphism. Now we treat cases. If $X$ has an infinite coarse component then Lemma \ref{lem:zandi}.(1) implies the map $K^{top}(X\times X) \rightarrow K^{top}(G(X))$ is the zero map. Now assume that $x \in K^{top}(G(X))$ maps to $0$ in $K_{0}(C^{*}_{r}(G(X))$. Then it maps to $0$ in $K^{top}(G(X)|_{\partial\beta X})$ as $\mu^{bdry}$ is an isomorphism. As the second line is exact, this implies it comes an element in $K^{top}(X\times X)$. As the map labelled $1$ is the zero map, $x$ must be $0$.

If $X$ has only finite coarse components, then by Lemma \ref{lem:zandi} the map $2$ is injective and so we can conclude injectivity of $\mu$ by the Five Lemma.
\end{proof}

We can also describe the obstructions to $\mu_{*}$ being an isomorphism when $X$ admits a fibred coarse embedding into Hilbert space.

\begin{proposition}\label{thm:MT2}
Let $X$ be a uniformly discrete metric space with bounded geometry such that $X$ admits a fibred coarse embedding into Hilbert space. Then the inclusion of $\mathcal{K}$ into $I_{G}$ induces an isomorphism on K-theory if and only if $\mu_{*}$ is an isomorphism. In addition, if $X$ has only finite coarse components then every ghost projection in $C^{*}X$ is compact if and only if $\mu_{0}$ is an isomorphism.
\end{proposition}
\begin{proof}
We consider the diagram from the proof of Theorem \ref{thm:mcor1}
$$
\xymatrix@=1em{
\ar[r] & K_{1}(C^{*}_{r}(G(X)|_{\partial\beta X}) \ar[r] & K_{0}(I_{G}) \ar[r]& K_{0}(C^{*}_{r}(G(X))) \ar[r]& K_{0}(C^{*}_{r}(G(X)|_{\partial\beta X})\ar[r] & K_{1}(I_{G}) \ar[r] & \\
\ar[r] & K_{1}^{top}(G(X)|_{\partial\beta X}) \ar[r] \ar[u]^{\mu^{bdry}}& K^{top}(X \times X) \ar[r]_{1}\ar[u]_{2}& K_{0}^{top}(G(X)) \ar[r]\ar[u]^{\mu}& K_{0}^{top}(G(X)|_{\partial\beta X}) \ar[r]\ar[u]^{\mu^{bdry}}& 0\ar[u]^{0} \ar[r] &
}
$$
When $X$ fibred coarsely embeds into Hilbert space, we know that $\mu^{bdry}$ is an isomorphism. The result then follows from the Five Lemma.

Now assume $X$ has only finite coarse components: As in Lemma \ref{lem:zandi} we will work using the generators of the metric coarse structure by considering the algebra $C^{*}_{R}X:=C^{*}_{\Delta_{R}}(X)$ and its intersection, $\mathcal{K}_{R}$, with the ideal of compact operators on $\ell^{2}X$: $\bigoplus_{i \in \Lambda_{R}}C^{*}(X_{i})$. We now want to construct a tracelike map, similar to the maps defined in \cite{higsonpreprint, explg1}, on each scale $R>0$.

For any $T \in \mathbb{C}[X]$ of propagation $R$, as in the proof of Lemma \ref{lem:zandi}.(2), we have a natural decomposition $T= \prod_{i \in \Lambda_{R}}T_{i}$, where each $T_{i} \in C^{*}X_{i}$ is the component of $T$ restricted to each of $X_{i}$ (Recall $X_{i}$ denotes the $i$-th component of the $0$-skeleton of the Rips complex $P_{E}X$). We define a map from $C^{*}_{R}X$ to $\frac{\prod_{i}C^{*}X_{i}}{\oplus C^{*}X_{i}}$ that sends each $T$ to $\overline{\prod T_{i}}$ and we denote the map it induces on K-theory by $Tr_{R}$.
\begin{eqnarray*}
Tr_{R}:K_{0}(C^{*}_{R}X) \rightarrow \frac{\prod_{i}K_{0}(C^{*}X_{i})}{\oplus_{i}K_{0}(C^{*}X_{i})}.
\end{eqnarray*}
To complete the proof, let $p$ be a noncompact ghost projection in $C^{*}X$, let $\lbrace T_{n} \rbrace_{n}$ be operators of propagation $R_{n}$ that approximate $p$ and let $\epsilon < \frac{1}{4}$. Observe that for some large $n$, there are quasiprojections $T_{n}$ such that $\Vert T_{n}-p \Vert \leq \epsilon$ that determine a K-theory class for $C^{*}_{\Delta_{R_{n}}}X$, which maps to a non-zero class under $Tr_{R_{n}}$. Suppose that $p$ was equivalent to some compact operator $q$ then we take representatives $q_{n} \in \mathcal{K}_{R_{n}}$ with $\Vert q_{n} - q \Vert \leq \epsilon$. Now, $Tr_{R_{n}}([q])=0$ but $Tr_{R_{n}}([T_{n}])\not = 0$: this gives the desired contradiction. 
 
If $\mu$ is an isomorphism, it follows that every ghost projection is equivalent to a compact operator on K-theory, that is $K_{*}(\mathcal{K}) \cong K_{*}(I_{G})$. Hence, any ghost projection in $C^{*}(X)$ vanishes under $Tr_{R}$ for a cofinal subset of $R>0$. This happens if and only if the ghost projection is compact \cite{explg1}.
\end{proof}

A natural corollary of this concerns coarsely embeddable spaces, where we can now show:

\begin{corollary}\label{thm:IT3}
If $X$ coarsely embeds into Hilbert space then $K_{*}(I_{G}) \cong K_{*}(\mathcal{K})$. \qed
\end{corollary}

\begin{remark}
This is motivated by some remarks of Roe \cite[Chapter 11]{MR2007488}, where he asks in what way coarse embeddability interacts with the existence of ghost operators. This was settled by Roe and Willett in \cite{ghostbusting}, where it was shown that the existence of non-compact ghost operators is equivalent to not having property A. Corollary \ref{thm:IT3} shows however that for spaces that do not have property A, but are coarsely embeddable into Hilbert space (such as spaces constructed in \cite{MR2899681,MR2920843}) that the ideals $\mathcal{K}$ and $I_{G}$ are not equal, but do have the same K-theory. For the spaces constructed in \cite{MR2899681,MR2920843} we can go further and say that whilst $I_{G}\not =\mathcal{K}$, there are no non-compact ghost projections in $I_{G}$ (This follows from Proposition \ref{thm:MT2}, as the spaces referenced above have only finite coarse components).
\end{remark}

\subsection{Box spaces of residually finite discrete groups}
Let $\Gamma$ be a finitely generated residually finite discrete group, and let $\mathcal{N}:=\lbrace N_{i} \rbrace_{i\in\mathbb{N}}$ be a family of nested finite index normal subgroups with trivial intersection. Fix a generating set $S$ for $\Gamma$. Then we can construct a metric space $\square \Gamma$, called the \textit{box space} of $\Gamma$ with respect to the family $\mathcal{N}$ by considering the coarse disjoint union of the sequence: $\lbrace Cay(\frac{\Gamma}{N_{i}},S) \rbrace_{i}$. 

It is well known \cite[Proposition 11.26]{MR2007488} that a coarse embedding of $\square \Gamma$ into Hilbert space implies that $\Gamma$ is a-T-menable. Using Theorem 2.2 from \cite{MR3116568}, it is possible to show that if $\Gamma$ is a-T-menable, then any box space of $\Gamma$ admits a fibred coarse embedding. The following Lemma will allow us to prove the converse of Theorem 2.2 from \cite{MR3116568} using Corollary \ref{Thm:MT1}.

\begin{lemma}\label{lem:cor2}
Let $\square \Gamma$ be a box space of a finitely generated residually finite group $\Gamma$. Then the boundary algebra $C(\partial\beta X)$ admits a $\Gamma$-invariant state.
\end{lemma}
\begin{proof}
Let $\omega \in \partial\beta \mathbb{N}$ and consider the function:
\begin{equation*}
\mu(f)=\lim_{\omega} \frac{1}{\vert X_{i} \vert}\sum_{x \in X_{i}}f(x_{i}).
\end{equation*}
Where the limit above is the ultralimit in $\omega$. Clearly, $\mu(1_{\beta X})=1$, $\mu$ is linear and positive, whence $\mu$ is a state on $C(\partial\beta X)$. Let $g \in \Gamma$, now we check invariance:
\begin{equation*}
\mu(g\circ f)=\lim \frac{1}{\vert X_{i} \vert}\sum_{x \in X_{i}}f(gx_{i}).
\end{equation*}
After relabelling the elements $x \in X_{i}$ by $g^{-1}x^{'}$, we now see that $\mu(f)=\mu(g \circ f)$. The result now follows.
\end{proof}

It is well known that the boundary groupoid associated to a box space $\square\Gamma$ decomposes as $\partial\beta \square \Gamma \rtimes \Gamma$, for a proof see \cite{mypub1}. We now recall Corollary 5.12 of \cite{MR3138486}:

\begin{proposition}\label{prop:cor}
Let $\Gamma$ be a discrete group acting on a space $X$ with an invariant probablity measure $\mu$. Then the action is a-T-menable if and only if $\Gamma$ is a-T-menable.\qed
\end{proposition}

Now Lemma \ref{lem:cor2} plus the above Proposition prove:

\begin{theorem}\label{thm:cor2}
Let $\Gamma$ be a residually finite discrete group. If the box space $\square \Gamma$ admits a fibred coarse embedding into Hilbert space then the group $\Gamma$ is a-T-menable.\qed
\end{theorem}

The author would like to remark that that Theorem \ref{thm:cor2} was proved independently using different methods in \cite{MR3105001}.


\bibliography{ref}
\end{document}